\documentclass[a4paper,11pt]{article}
\usepackage[latin1]{inputenc}
\usepackage{graphicx}
\usepackage{bm}
\usepackage{epsfig}
\usepackage{epsfig}
\usepackage{amsmath}
\usepackage{amsfonts}
\usepackage{amssymb, longtable}
\usepackage{relsize}

\newtheorem{Thm}{Theorem}
\newtheorem{Prop}[Thm]{Proposition}
\newtheorem{Def}[Thm]{Definition}

\usepackage{color}
\begin{document}

\begin{center}
{\LARGE\bf An extension of Azzalini's method}\\
\ \\ \ \\ \ \\ \ \\
{\large\bf Filippo Domma}\\
\vspace{0.5 cm}
{
Department of Economics, Statistics and Finance,\\
University of Calabria, Italy\\ [1cm]}

{\large\bf Bo\v zidar V. Popovi\'c}\\
\vspace{0.5 cm}
{
University of Montenegro, Faculty of Science and Mathematics,\\
Podgorica, Montenegro \\ [1cm]}

{\large\bf Saralees Nadarajah}\\
\vspace{0.5 cm}
{
School of Mathematics, University of Manchester,\\
Manchester, UK\\ [1cm]}
\end{center}

\begin{abstract}
The aim of this paper is to extend Azzalini's method.
This extension is done in two stages:
consider two dependent and non-identically distributed random variables say $X_1$ and $X_2$;
model the dependence between $X_1$ and $X_2$ by a copula.
To illustrate the new method, we assume $X_1$ and $X_2$ are exponential random variables.
This assumption leads to a  new distribution called the
Generalized Weighted Exponential Distribution (GWED),
a generalization of Gupta and Kundu (2009)'s Weighted Exponential Distribution (WED).
Some mathematical properties of the GWED are derived, and its parameters estimated by maximum likelihood.
The GWED is applied to biochemical data sets showing its good performance compared to the WED.
\end{abstract}

\vspace{0.4cm}


\vspace{0.4cm}


\section{Introduction}

In real life, there are many data sets that are  asymmetric, multimodal and heavy tailed.
This has motivated  many researchers to develop non-normal and/or skewed distributions.
In the statistical literature, there are various techniques to build non-normal and/or skewed distributions.
Nowadays, the most widely used technique for introducing asymmetry in a symmetric distribution
is that due to Azzalini (1985).
With reference to the normal distribution, this technique can be described as follows:
a random variable $Z$ is said to have the skew-normal distribution with parameter $\alpha$,
written as $Z \sim SN(\alpha)$, if its probability density function (pdf) is
$f(z;\alpha)=2\Phi(\alpha z) \phi (z)$,
where $\Phi(\cdot)$ denotes the standard normal cumulative distribution function (cdf),
$\phi(\cdot)$ denotes the standard normal pdf,
and $z$ and $\alpha$ are real numbers (Azzalini, 1985).
An enormous literature exists on the study of the skew-normal distribution
and its extension to the multivariate case.

Azzalini's method can also be described in terms of conditional distributions:
let $X_i$, $i = 1, 2$ be independent random
variables with pdfs $f_i\left(x_i\right)$ and cdfs $F_i\left(x_i\right)$.
Then, the conditional pdf of $X=X_{1}$ given $ \alpha X_{1} > X_2$ is
\begin{eqnarray}
\displaystyle
f_{X}(x)=\frac {\displaystyle f_{1}(x) F_{2}(\alpha x)}{\displaystyle P \left(\alpha X_{1} > X_{2} \right)}.
\label{WAzz}
\end{eqnarray}
Observe that $P \left(\alpha X_{1} > X_{2} \right) = \displaystyle \int^{+\infty}_{0}
\left[\int^{\alpha x_1}_{0}f_2 \left(x_2\right) {\rm d}x_2 \right]f_1 \left(x_1\right) {\rm d}x_1 = E_{X_{1}}
\left[ F_{2} (\alpha X) \right]$, where $E_{X_{1}}(\cdot)$ denotes the
expectation with respect to $X_1$.
Equation (\ref{WAzz}) can be interpreted as a weighted distribution
with weight function $w(x) = F_{2}(\alpha x)$.
If we set $f_1\left(x_1\right) = f_2\left(x_2\right)=\phi(x)$, i.e., $X_1$ and $X_2$ are standard normal random variables,
we obtain the skew-normal distribution.

In the literature, (\ref{WAzz}) has been used to construct new skewed distributions from a given symmetric distribution,
for example, skew-$t$ (Gupta et al., 2002), skew-Cauchy (Arnold and Beaver, 2000b; Gupta et al., 2002),
skew-Laplace (Gupta et al., 2002; Aryal and Nadarajah, 2005)
and skew-logistic (Gupta et al., 2002; Nadarajah, 2009).
However, there is little work on the use of Azzalini's method for non-symmetric distributions.
Gupta and Kundu (2009) took $X_1$ and $X_2$ in (\ref{WAzz})
to be independent and identical exponential random variables with scale parameter $\lambda$, giving
\begin{eqnarray*}
\displaystyle
f_{X}(x)=\frac {\displaystyle \alpha +1}{\displaystyle \alpha} \lambda e^{-\lambda x} \left(1-e^{-\alpha\lambda x}\right).
\end{eqnarray*}
Gupta and Kundu (2009) called this the Weighted Exponential Distribution (WED).
Shakhatreh (2012) studied a two-parameter version of the WED.
Mahdy (2011, 2013) proposed weighted gamma and weighted Weibull distributions.

The aim of this paper is to extend Azzalini's method in two stages:
take $X_1$ and $X_2$ to be dependent and non-identically distributed random variables;
model their dependence using a copula.
After a general discussion on the potentials,
we illustrate this method by assuming $X_1$ and $X_2$ are exponential random variables.
This assumption leads us to a  new distribution called the Generalized Weighted Exponential Distribution (GWED), a generalization of the WED.
Although the GWED is defined using Azzalini's method in terms of weighted distribution,
we will show that it can also be interpreted as a hidden truncated distribution (Arnold and Beaver, 2000a).
Moreover, the GWED can be interpreted as a finite mixture.
The mixture representation enables us to derive mathematical properties of the GWED like
its cdf, moments, and the moment generating function.
The skewness of the WED due to Gupta and Kundu (2009) takes values in $\left[\sqrt{2},2\right]$
while that of the two-parameter WED due to Shakhatreh (2012) takes values in $\left[ \frac {2}{\sqrt{3}}, 2 \right]$.
The GWED allows for wider values for skewness.
Another importance feature is that the GWED allows for non-monotonic hazard rate functions (hrfs).
The hrfs of the WED are always monotonic.

The paper is organized as follows.
The extension of Azzalini's method is described in Section 2.
Details (including mathematical properties, estimation issues and applications) of an important special case are given in Section 3.

\section{The extension of Azzalini's method}

Unlike Azzalini's method, we consider two dependent and non-identically distributed
random variables.
This extension can be used to construct any distribution.

First note that the denominator in (\ref{WAzz}) in the case of dependence of $X_1$ and $X_2$ can be expressed as
\begin{eqnarray}
\displaystyle
P\left(\alpha X_{1}>X_{2}\right)
&=&
\displaystyle
\int^{+\infty}_{0} \int^{\alpha x_1}_{0} f \left(x_{1}, x_{2}\right) {\rm d}x_{1} {\rm d}x_{2}
\nonumber
\\
&=&
\displaystyle
\int^{+\infty}_{0} \left\{ \int^{\alpha x_1}_{0}  c\left( F_{1} \left( x_{1} \right), F_{2} \left( x_{2} \right) \right)
f_{2} \left(x_{2}\right) {\rm d}x_{2} \right\} f_{1} \left(x_{1}\right) {\rm d}x_{1}
\nonumber
\\
&=&
\displaystyle
E_{X_{1}}\left[ \int^{\alpha X_1}_{0}  c\left( F_{1} \left( X_{1} \right), F_{2} \left(x_{2} \right) \right) f_{2} \left(x_{2}\right) {\rm d}x_{2} \right],
\label{ssm}
\end{eqnarray}
where $c\left( F_{1} \left( x_{1} \right), F_{2} \left(x_{2} \right) \right) =
\frac {\partial^{2} C\left( F_{1} \left( x_{1} \right), F_{2} \left( x_{2} \right) \right)}
{\partial F_{1} \left( x_{1} \right) \partial F_{2} \left( x_{2} \right)}$ is the copula pdf.
Accordingly, the weight function is $w(x) = \displaystyle\int^{\alpha x_1}_{0}  c\left( F_{1} \left(  x_{1} \right), F_{2} \left( x_{2} \right) \right)
f_{2} \left( x_{2} \right) {\rm d}x_{2}$.
If $\alpha=1$, (\ref{ssm}) reduces to the stress-strength model
widely studied in the literature, see, for example, Kotz et al. (2003);
using a copula, Domma and Giordano (2013) have recently highlighted the role of dependence between stress and strength
on a reliability measure defined by (\ref{ssm}).

We are now able to provide

\begin{Def}
Let $X_i$, $i = 1, 2$ be dependent random variables with pdfs
$f_{i} \left( x_{i};\xi_{i} \right)$, cdfs $F_{i} \left(x_{i};\xi_{i}\right)$
and joint pdf $f \left( x_{1},x_{2}; \xi_{1}, \xi_{2}, \theta \right) =
c \left( F_{1} \left( x_{1}; \xi_{1} \right), F_{2} \left( x_{2};\xi_{2} \right); \theta \right)
f_{1} \left( x_{1};\xi_{1} \right) f_{2} \left( x_{2};\xi_{2} \right)$,
where $c(\cdot, \cdot)$ is a copula pdf.
Then, the random variable $X=X_{1} | \alpha X_{1}>X_{2}$
is said to have a Generalized Weighted Distribution if its pdf is
\begin{eqnarray}
\displaystyle
f^{w} \left( x; \xi_{1}, \xi_{2}, \alpha, \theta \right) =
\frac {\displaystyle f_{1} \left( x; \xi_{1} \right) \int_{-\infty}^{\alpha x} c \left( F_{1} \left( x;\xi_{1} \right),
F_{2} \left( x_{2}; \xi_{2} \right); \theta \right) f_{2} \left( x_{2}; \xi_{2} \right) {\rm d}x_{2}}
{\displaystyle E_{X_{1}}\left[ \int^{\alpha X_1}_{-\infty}  c\left( F_{1} \left( X_{1} \right), F_{2} \left( x_{2} \right); \xi_{1} \right)
f_{2} \left( x_{2}; \xi_{2} \right) {\rm d}x_{2} \right]}.
\label{newpdf}
\end{eqnarray}
\end{Def}

The use of copula is motivated by the fact that it allows for
the dependence structure to be treated separately from the marginal components of a joint distribution.
Various forms of dependence structures (linear, non-linear, tail dependence, etc) can be used.
Also $X_{1}$ and $X_{2}$ need not necessarily belong to the same parametric family.

To better understand the role of (\ref{newpdf}) in defining new distributions,
we now consider a recent expansion applicable to a wide range of bivariate copulas (Nadarajah, 2014):
\begin{eqnarray}
\displaystyle
C\left( F_{1} \left( x_{1}; \xi_{1} \right), F_{2} \left( x_{2}; \xi_{2} \right); \theta  \right) =
\sum^{p}_{i=1} \gamma_{i} \left[F_{1} \left( x_{1}; \xi_{1} \right) \right]^{a_{i}} \left[F_{2} \left( x_{2}; \xi_{2} \right)\right]^{b_{i}},
\label{expC}
\end{eqnarray}
where $p\geq 1$ and $\left\{ \left( \gamma_{i}, a_{i}, b_{i} \right) : i\geq 1 \right\}$ are real numbers.
The corresponding joint pdf is
\begin{eqnarray*}
\displaystyle
c\left( F_{1} \left( x_{1}; \xi_{1} \right), F_{2} \left(  x_{2}; \xi_{2} \right); \theta  \right) =
\sum^{p}_{i=1} \gamma_{i} a_{i}b_{i}\left[F_{1} \left( x_{1}; \xi_{1} \right)\right]^{a_{i}-1}
\left[F_{2} \left( x_{2}; \xi_{2} \right) \right]^{b_{i}-1}.
\end{eqnarray*}
Then, the denominator in (\ref{newpdf}) reduces to
\begin{eqnarray*}
\displaystyle
P\left( \alpha X_{1}>X_{2} \right) = \sum^{p}_{i=1} \gamma_{i} a_{i}
E_{X_{1}}\left[ \left( F_{1} \left( X_{1}; \xi_{1} \right) \right)^{a_{i}-1}
\left( F_{2} \left(\alpha X_{1}; \xi_{2} \right) \right)^{b_{i}}   \right]\,,
\end{eqnarray*}
while weight function is
\begin{eqnarray*}
\displaystyle
w(x)= \sum^{p}_{i=1} \gamma_{i} a_{i} \left(F_{1} \left( x; \xi_{1} \right)\right)^{a_{i}-1}
\left( F_{2} \left( \alpha x; \xi_{2} \right) \right)^{b_{i}}.
\end{eqnarray*}
Consequently, the new pdf is
\begin{eqnarray}
\displaystyle
f^{w} \left( x;\xi_{1}, \xi_{2}, \alpha, \theta \right) =
\overline{K} f_{1} \left( x; \xi_{1} \right)
\sum^{p}_{i=1} \gamma_{i} a_{i} \left(F_{1} \left( x; \xi_{1} \right) \right)^{a_{i} - 1}
\left( F_{2} \left( \alpha x; \xi_{2} \right) \right)^{b_{i}},
\label{pdfexpC}
\end{eqnarray}
where
\begin{eqnarray*}
\displaystyle
\overline{K}=\left\{\sum^{p}_{i=1} \gamma_{i} a_{i}
E_{X_{1}}\left[ \left( F_{1} \left( X_{1}; \xi_{1} \right) \right)^{a_{i}-1}
\left( F_{2} \left( \alpha X_{1}; \xi_{2} \right) \right)^{b_{i}}   \right]\right\}^{-1}.
\end{eqnarray*}
The properties of the new pdf (\ref{pdfexpC}) depend
on the specified forms of the cdfs $F_{1} \left( x_{1}; \xi_{1} \right)$ and $F_{2} \left( x_{2}; \xi_{2} \right)$.
In the next section, we study a simple but an important special case that generalizes a distribution proposed by Gupta and Kundu (2009).

\section{A special case: Generalized weighted exponentiated distribution}

Here, we illustrate an application of the extension of Azzalini's method.
We take $X_i$, $i = 1, 2$ to be exponential random variables with cdfs $F_i (x) = \lambda_i e^{-\lambda_i x},$
and $C$ to be the Farlie-Gumbel-Morgenstern copula (Morgenstern, 1956) defined by
\begin{eqnarray*}
C \left( u_1, u_2 \right) =  u_{1} u_{2} \left\{ 1+\theta \left[1- u_{1}\right]\left[1-u_{2}\right]  \right\}.
\end{eqnarray*}
These choices yield a new distribution called the Generalized Weighted Exponential Distribution (GWED).
This distribution generalizes the WED due to Gupta and Kundu (2009).
We will show that the GWED can be  interpreted as a hidden truncated distribution (Arnold and Beaver, 2000a)
and as a finite mixture.
We will also discuss the behavior of the pdf and the hrf of the GWED.

The Farlie-Gumbel-Morgenstern copula is a special case of (\ref{expC}) for
$p=4$, $\gamma_{1}=1+\theta$, $a_{1}=1$, $b_{1}=1$, $\gamma_{2}=-\theta$, $a_{2}=1$, $b_{2}=2$,
$\gamma_{3}=-\theta$, $a_{3}=2$, $b_{3}=1$, $\gamma_{4}=\theta$, $a_{4}=2$ and $b_{4}=2$.
Consequently, after simple algebra, we find that
$\overline{K}=\left\{ \frac {\lambda_2 \alpha}{\lambda_1+\lambda_2 \alpha} +
\theta \lambda_1 \left[\frac {2}{2\lambda_1+\lambda_2 \alpha} +
\frac {1}{\lambda_1+2\lambda_2 \alpha}-\frac {2}{\lambda_1+\lambda_2 \alpha}\right] \right\}^{-1}$.
Now, we can provide the definition of the GWED.

\begin{Def}
Let $X_i$, $i = 1, 2$ be exponential random variables with pdfs
$f_{i} \left( x_{i};\lambda_{i} \right) = \lambda_{i} e^{-\lambda_{i} x_{i}}$,
cdfs $F_{i} \left( x_{i};\lambda_{i} \right) = 1 - e^{-\lambda_{i} x_{i}}$
and the joint pdf
\begin{eqnarray*}
\displaystyle
f \left( x_{1}, x_{2}; \lambda_{1}, \lambda_{2}, \theta \right) =
\left\{1+\theta\left[1-2 e^{-\lambda_{1} x_{1}}\right] \left[1-2 e^{-\lambda_{2} x_{2}}\right]\right\}
f_{1} \left( x_{1}; \lambda_{1} \right) f_{2} \left( x_{2}; \lambda_{2} \right).
\end{eqnarray*}
Then, the random variable $X=X_{1} | \alpha X_{1}>X_{2}$ is said to have the GWED if its pdf is
\begin{eqnarray}
\label{pdfGWE}
\displaystyle
f^{w} \left( x;\lambda_{1}, \lambda_{2}, \alpha, \theta \right) =
K e^{-\lambda_1 x}\left(1-e^{-\lambda_2 \alpha x}\right)
\left[1-\theta e^{-\lambda_2 \alpha x} \left(1-2e^{-\lambda_1 x}\right)\right],
\end{eqnarray}
where $K=\lambda_{1} \overline{K}$.
\end{Def}

It is easy to see that the WED due to
Gupta and Kundu (2009)
is the special case of the GWED for  $\lambda_1=\lambda_2=\lambda$ and $\theta=0$.

Proposition 3 shows that $f^{w} \left(x;\lambda_{1}, \lambda_{2}, \alpha, \theta \right)$
in (\ref{pdfGWE}) can be interpreted as a hidden truncated pdf (Arnold and Beaver, 2000a).

\begin{Prop}
Suppose $Z$ and $Y$ are two dependent positive random variables with joint pdf
\begin{eqnarray*}
\displaystyle
f_{Z,Y} \left( z, y; \lambda_{1}, \lambda_{2}, \theta \right) =
\lambda_{1} \lambda_{2} z e^{-\left( \lambda_{1}+\lambda_{2}y \right) z}
\left\{1+\theta \left[ 1-2 e^{-\lambda_{1} z} \right] \left[ 1-2 e^{-\lambda_{2} y z} \right] \right\}.
\end{eqnarray*}
Then, the pdf of $X=Z | Y\leq \alpha$ is equal to $f^{w} \left( x;\lambda_{1}, \lambda_{2}, \alpha, \theta \right)$ in (\ref{pdfGWE}).
\end{Prop}

\noindent
\textbf{Proof.}
First, observe that the conditional cdf of $X=Z | Y\leq \alpha$ is
\begin{eqnarray*}
\displaystyle
F_{X=Z|Y<\alpha}(x)=\frac {\displaystyle \int^{x}_{0}  \int^{\alpha}_{0} f_{Z,Y}(z,y) {\rm d}y {\rm d}z}
{\displaystyle \int^{\alpha}_{0} f_{Y}(y) {\rm d}y}.
\end{eqnarray*}
After simple algebra, we have
\begin{eqnarray*}
\displaystyle
\int^{\alpha}_{0}f_{Z,Y} \left(z,y; \lambda_{1}, \lambda_{2}, \theta\right) {\rm d}y
&=&
\displaystyle
\lambda_{1}\left\{ (1+\theta) e^{-\lambda_{1} z} \left(1-e^{-\lambda_{2} \alpha z} \right) -
\theta e^{-\lambda_{1} z} \left(1-e^{-2\lambda_{2} \alpha z} \right) \right.
\\
&&
\displaystyle
\left.- 2\theta e^{-2 \lambda_{1} z} \left(1-e^{-\lambda_{2} \alpha z} \right) +
2\theta e^{-2 \lambda_{1} z} \left(1-e^{-2\lambda_{2} \alpha z} \right) \right\}
\end{eqnarray*}
and
\begin{eqnarray}
&&
\displaystyle
\int^{x}_{0}  \int^{\alpha}_{0}f_{Z,Y} \left( z,y; \lambda_{1}, \lambda_{2}, \theta \right) {\rm d}y =
\lambda_{1}  \left\{
(1+\theta) \left[ \frac {\displaystyle 1-e^{-\lambda_{1} x}}{\displaystyle \lambda_{1}} -
\frac {\displaystyle 1-e^{-\left( \lambda_{1}+\lambda_{2} \alpha \right)x}}{\displaystyle \lambda_{1}+\lambda_{2} \alpha}  \right]
\right.
\nonumber
\\
&&
\displaystyle
-\left. \theta \left[  \frac {\displaystyle 1 - e^{-\lambda_{1} x}}{\displaystyle \lambda_{1}} -
\frac {\displaystyle 1-e^{-\left( \lambda_{1}+2\lambda_{2} \alpha \right)x}}{\displaystyle \lambda_{1}+ 2\lambda_{2} \alpha} \right] -
2\theta \left[
\frac {\displaystyle 1-e^{-2\lambda_{1} x}}{\displaystyle 2\lambda_{1}} -
\frac {\displaystyle 1-e^{-\left( 2\lambda_{1}+\lambda_{2} \alpha \right)x}}{\displaystyle 2\lambda_{1}+ \lambda_{2} \alpha} \right]
\right.
\nonumber
\\
&&
\displaystyle
+\left. 2\theta \left[ \frac {\displaystyle 1-e^{-2\lambda_{1} x}}{\displaystyle 2\lambda_{1}} -
\frac {\displaystyle 1-e^{-2 \left( \lambda_{1}+\lambda_{2} \alpha \right) x}}
{\displaystyle 2\left( \lambda_{1}+ \lambda_{2} \alpha \right)} \right] \right\}.
\label{num}
\end{eqnarray}
Moreover, the marginal pdf of $Y$ is
\begin{eqnarray*}
\displaystyle
f_{Y} \left( y; \lambda_{1}, \lambda_{2}, \theta \right) =
\lambda_{1} \lambda_{2} \left\{ \frac {\displaystyle 1+2\theta}{\displaystyle \left( \lambda_{1}+\lambda_{2} y \right)^2} -
\frac {\displaystyle 2\theta}{\displaystyle \left( \lambda_{1}+2\lambda_{2} y \right)^2} -
\frac {\displaystyle 2\theta}{\displaystyle \left( 2\lambda_{1}+\lambda_{2} y \right)^2} \right\}.
\end{eqnarray*}
Consequently,
\begin{eqnarray}
\displaystyle
P\left(Y\leq \alpha\right) = \frac {\displaystyle \lambda_{2} \alpha}{\displaystyle \lambda_{1}+\lambda_{2} \alpha} +
\lambda_{1} \theta
\left[ \frac {\displaystyle 2}{\displaystyle 2\lambda_{1}+\lambda_{2} \alpha} +
\frac {\displaystyle 1}{\displaystyle \lambda_{1}+2\lambda_{2} \alpha} -
\frac {\displaystyle 2}{\displaystyle \lambda_{1}+\lambda_{2} \alpha} \right].
\label{den}
\end{eqnarray}
By combining (\ref{num}) with (\ref{den}), we obtain the cdf of $X|Y\leq \alpha$,
i.e., $F_{X=Z|Y\leq \alpha} \left( x;\lambda_{1}, \lambda_{2}, \alpha, \theta \right)$.
By differentiating with respect to $x$, it  can be verified
that the conditional pdf $f_{X=Z|Y\leq \alpha} \left( x;\lambda_{1}, \lambda_{2}, \alpha, \theta \right)$ is equal to (\ref{pdfGWE}).
The proof is complete.
\
$\Box$

\begin{figure}
\includegraphics[width=0.75\textwidth]{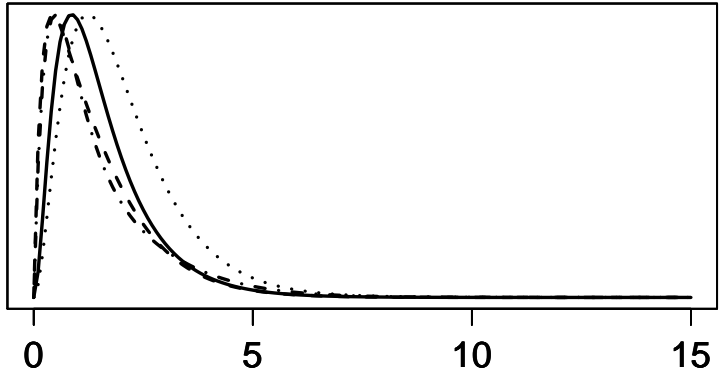}
\caption{Examples of the GWED pdf for
$\left(\lambda_1,\lambda_2,\alpha,\theta\right)=(1,1,1,-0.99)$ (solid line),
$\left(\lambda_1,\lambda_2,\alpha,\theta\right)=(1,1,1,0.99)$ (dashed line),
$\left(\lambda_1,\lambda_2,\alpha,\theta\right)=(1,1,0.5,-0.99)$ (dotted line),
$\left(\lambda_1,\lambda_2,\alpha,\theta\right)=(1,1,0.5,0.99)$ (dot-dashed line).}
\label{pdfss}
\end{figure}

It is easy to see that
$\lim_{x\rightarrow 0^+}f^{w} \left( x;\lambda_{1}, \lambda_{2}, \alpha, \theta \right)=
\lim_{x\rightarrow +\infty}f^{w} \left( x;\lambda_{1}, \lambda_{2}, \alpha, \theta \right) = 0$.
Therefore, there exists at least one $x>0$, say $x_0$, such that
$\left. \frac {\partial f^{w} \left( x;\lambda_{1}, \lambda_{2}, \alpha, \theta \right)}{\partial x} \right|_{x=x_{0}} =0$.
Figure \ref{pdfss} plots \eqref{pdfGWE} for different parameter values.

Proposition 4 interprets (\ref{pdfGWE}) as a finite mixture of exponential pdfs $f_{j} (x; \tau) = \tau e^{-\tau x}$.

\begin{Prop}
The pdf $f^{w} \left( x; \lambda_{1}, \lambda_{2}, \alpha, \theta \right)$ can be expressed as
\begin{eqnarray}
\displaystyle
f^{w} \left( x;\lambda_{1}, \lambda_{2}, \alpha, \theta \right)
&=&
\displaystyle
p_{1} f_{1} \left( x; \lambda_{1} \right) + p_{2} f_{2} \left( x; \lambda_{1}+\lambda_{2}\alpha \right) +
p_{3} f_{3} \left( x; 2\lambda_{1}+\lambda_{2}\alpha \right)
\nonumber
\\
&&
\displaystyle
+p_{4} f_{4} \left( x; \lambda_{1}+2\lambda_{2}\alpha \right) +
p_{5} f_{5} \left( x; 2 \left( \lambda_{1}+\lambda_{2}\alpha \right) \right),
\label{pdfGWEmix}
\end{eqnarray}
where $p_{1} = \frac {K}{\lambda_{1}}$,
$p_{2}=-\frac {K(1+\theta)}{\lambda_{1}+\lambda_{2}\alpha}$,
$p_{3}=\frac {2K\theta}{2\lambda_{1}+\lambda_{2}\alpha}$,
$p_{4}=\frac {K\theta}{\lambda_{1}+2\lambda_{2}\alpha}$ and
$p_{5}=-\frac {K\theta}{\lambda_{1}+\lambda_{2}\alpha}$.
\end{Prop}

\noindent
\textbf{Proof.}
Follows by simple algebra.
\
$\Box$

For details about  mixtures, we refer readers to  Titterington et al. (1985, page 50).
A mixture of distributions  is a valid distribution if the sum of weights is equal to $1$.
In the case of (\ref{pdfGWEmix}), it is straightforward to verify that $\displaystyle \sum^{5}_{j=1}p_{j}=1$.
The mixture representation (\ref{pdfGWEmix}) enables us to derive mathematical properties
of the GWED like its cdf, moments, and moment generating function.
For example, the cdf of the GWED can be expressed as
\begin{eqnarray}
\displaystyle
F^{w} \left( x;\lambda_{1}, \lambda_{2}, \alpha, \theta \right)
&=&
p_{1} F_{1} \left(x;\lambda_{1}\right)+p_{2} F_{2}\left(x;\lambda_{1}+\lambda_{2}\alpha\right)+
p_{3} F_{3} \left(x;2\lambda_{1}+\lambda_{2}\alpha\right)
\nonumber
\\
&&
\displaystyle
+p_{4} F_{4} \left( x;\lambda_{1}+2\lambda_{2}\alpha \right) +
p_{5} F_{5} \left( x; 2 \left( \lambda_{1}+\lambda_{2}\alpha \right) \right).
\nonumber
\end{eqnarray}
An alternative expression for the cdf can be obtained from  \eqref{pdfGWEmix} by using the Taylor series for exponential functions:
\begin{eqnarray}
\label{P1}
\displaystyle
F^{w} \left( x;\lambda_{1}, \lambda_{2}, \alpha, \theta \right) =
\sum_{i\ge1}a_i x^i = \sum_{i\ge0}a_{i+1} x^{i+1}=x\sum_{i\ge0}d_i x^{i},
\end{eqnarray}
where
\begin{eqnarray*}
\displaystyle
d_i = a_{i+1}
&=&
\displaystyle
\frac {\displaystyle (-1)^{i+1}}{\displaystyle (i+1)!}
\Bigg[p_{1}\lambda_1^{i+1}+p_{2} \left( \lambda_1+\lambda_2\alpha \right)^{i+1}
+p_{3} \left(2\lambda_1+\lambda_2\alpha\right)^{i+1}
\\
&&
\qquad \qquad \qquad \qquad \qquad \qquad
\displaystyle
+p_{4} \left(\lambda_1+2\lambda_2\alpha\right)^{i+1}+
+2^{i+1} p_{5}  \left(\lambda_1+\lambda_2\alpha\right)^{i+1}\Bigg].
\end{eqnarray*}
The hrf  of the GWED  given by $h^{w} \left( x;\lambda_{1}, \lambda_{2}, \alpha, \theta \right) =
\frac {f^{w} \left( x;\lambda_{1}, \lambda_{2}, \alpha, \theta \right)}
{1 - F^{w} \left( x;\lambda_{1}, \lambda_{2}, \alpha, \theta \right)}$ can also be expressed using the mixture representation  (\ref{pdfGWEmix}).
It will take a complicated expression.
However, the behavior of the hrf can be easily assessed:
we can verify that $\lim_{x\rightarrow 0^+}h^{w} \left( x;\lambda_{1}, \lambda_{2}, \alpha, \theta \right) = 0$
and $\lim_{x\rightarrow +\infty}h^{w} \left( x;\lambda_{1}, \lambda_{2}, \alpha, \theta \right)=\lambda_1$,
where the latter follows by L'Hospital rule.

\begin{figure}[!htbp]
(a)\hspace{9cm}    (b)

\includegraphics[width=0.45\textwidth]{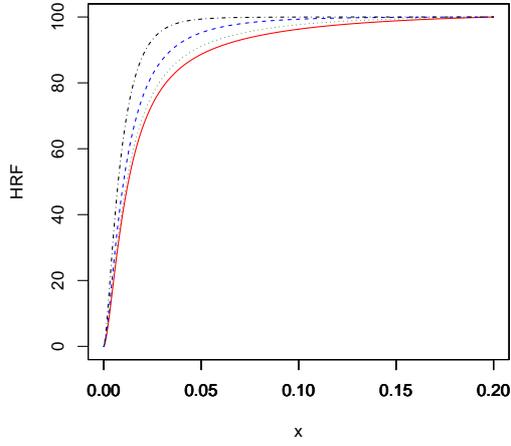}
\includegraphics[width=0.45\textwidth]{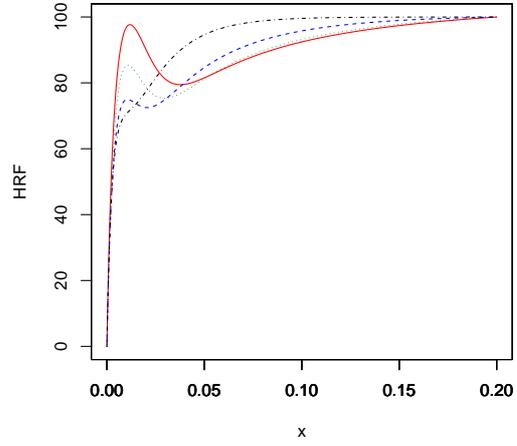}
\caption{(a) The hrf  of the GWED
for $\theta<0$, $\alpha=0.5$ (red line), $\theta<0$, $\alpha=1$ (green line),
$\theta<0$, $\alpha=2$ (blue line) and $\theta<0$, $\alpha=5$ (black line);
(b) The hrf  of the GWED for $\theta>0$, $\alpha=0.5$ (red line),
$\theta>0$, $\alpha=1$ (green line), $\theta>0$, $\alpha=2$ (blue line) and $\theta>0$, $\alpha=5$ (black line).}
\label{hrfss1a}
\end{figure}

\begin{figure}[!htbp]
(a) \hspace{9cm}  (b)

\includegraphics[width=0.45\textwidth]{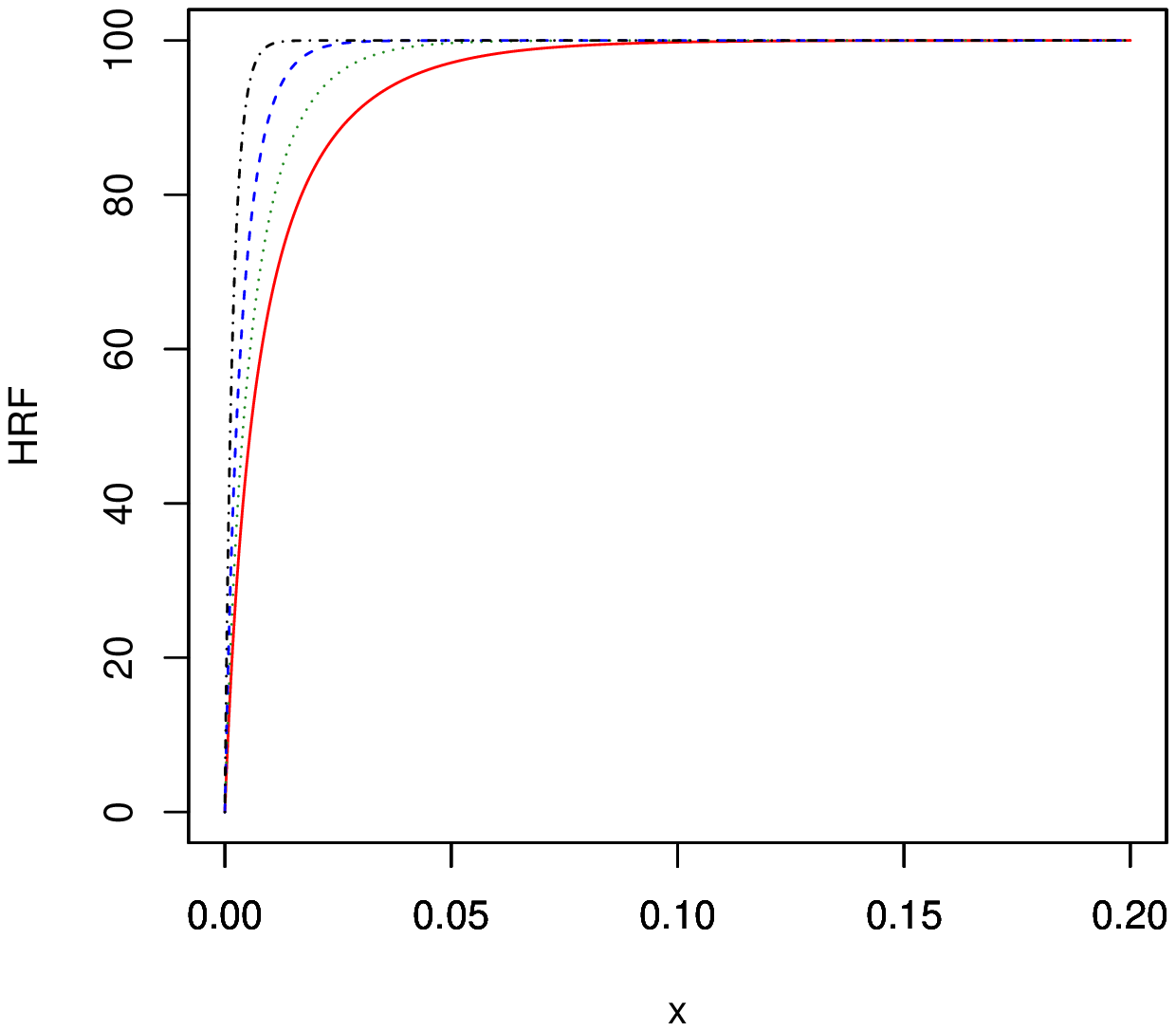}
\includegraphics[width=0.45\textwidth]{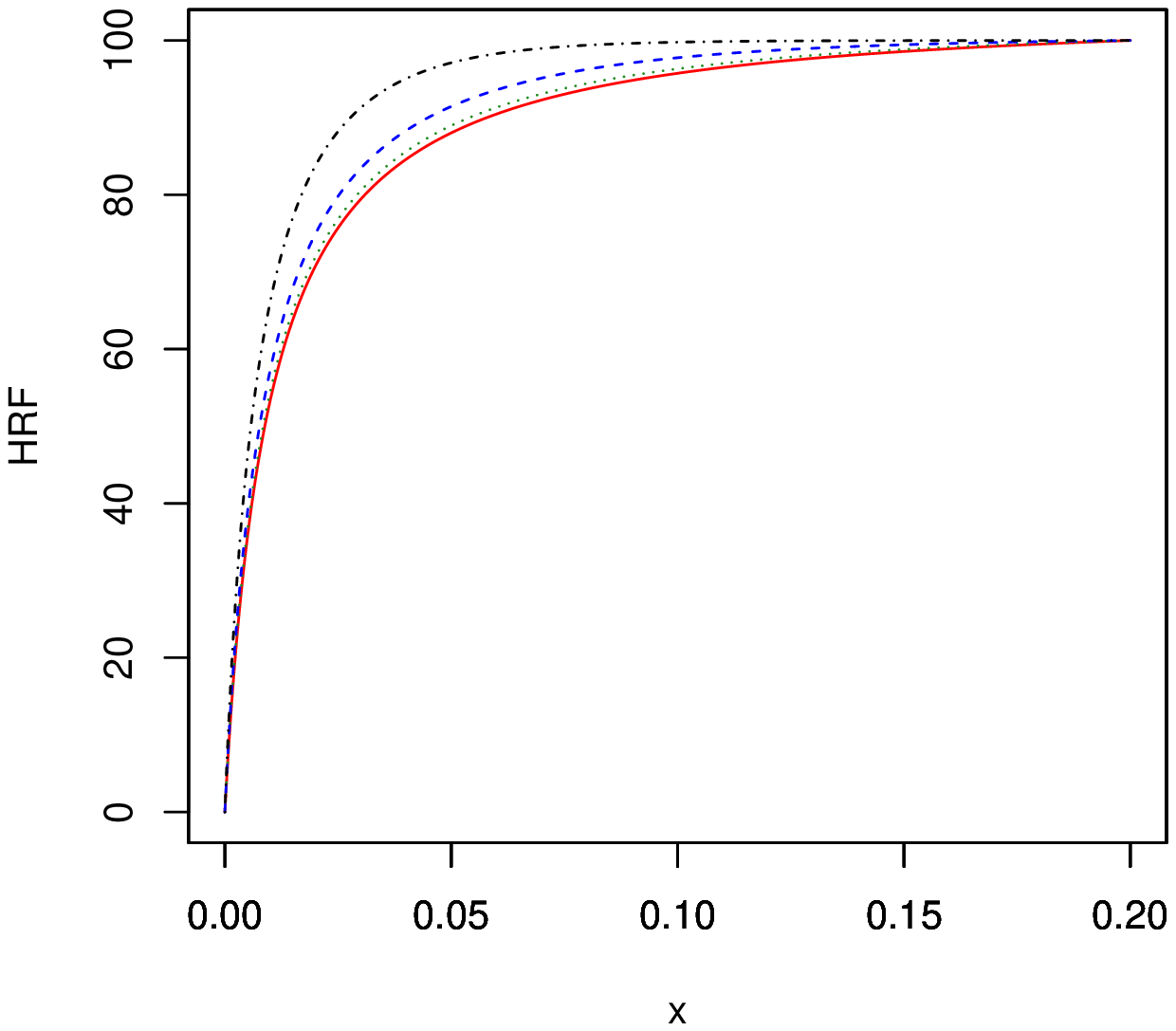}
\caption{(a) The hrf of the WED for $\alpha=0.5$ (red line),
$\alpha=1$ (green line), $\alpha=2$ (blue line) and $\alpha=5$ (black line);
(b)  The hrf of the GWED for $\theta=0$, $\alpha=0.5$ (red line), $\theta=0$, $\alpha=1$ (green line),
$\theta=0$, $\alpha=2$ (blue line) and $\theta=0$, $\alpha=5$ (black line).}
\label{hrfss2a}
\end{figure}

Figures \ref{hrfss1a} and \ref{hrfss2a} illustrate the behavior of the hrf for different parameter values.
We see that $h^{w} \left( x;\lambda_{1}, \lambda_{2}, \alpha, \theta \right)$ of the GWED
is more flexible than the hrf of the WED.
The hrfs of the GWED appear monotonic for negative dependence and non-monotonic for positive dependence.

\subsection{Some mathematical properties}

Let $X$ be a GWE random variable.
Proposition 5 derives the $n$th moment and the moment generating function of $X$.

\begin{Prop}
The $n$th  moment and the moment generating function of $X$ can be expressed as
\begin{eqnarray*}
\displaystyle
\mu_{n}^{\prime}=n! \left\{\frac {\displaystyle p_{1}}{\displaystyle \lambda_{1}^{n}} +
\frac {\displaystyle p_{2}}{\displaystyle \left( \lambda_{1}+\lambda_{2}\alpha \right)^n} +
\frac {\displaystyle p_{3}}{\displaystyle \left( 2\lambda_{1}+\lambda_{2}\alpha \right)^n} +
\frac {\displaystyle p_{4}}{\displaystyle \left( \lambda_{1}+2\lambda_{2}\alpha \right)^n} +
\frac {\displaystyle p_{5}}{\displaystyle 2^{n} \left( \lambda_{1}+\lambda_{2}\alpha \right)^n} \right\}
\end{eqnarray*}
and
\begin{eqnarray*}
\displaystyle
\mathcal{M}(t) =  \frac {\displaystyle p_{1}\lambda_{1}}{\displaystyle \lambda_{1} - t} +
\frac {\displaystyle p_{2} \left( \lambda_{1}+\lambda_{2}\alpha \right)}{\displaystyle \lambda_{1}+\lambda_{2}\alpha -t} +
\frac {\displaystyle p_{3} \left( 2\lambda_{1}+\lambda_{2}\alpha \right)}{\displaystyle 2\lambda_{1}+\lambda_{2}\alpha -t} +
\frac {\displaystyle p_{4} \left( \lambda_{1}+2\lambda_{2}\alpha \right)}{\displaystyle \lambda_{1}+2\lambda_{2}\alpha -t} +
\frac {\displaystyle 2p_{5} \left( \lambda_{1}+\lambda_{2}\alpha \right)}{\displaystyle 2\lambda_{1}+2\lambda_{2}\alpha -t},
\end{eqnarray*}
respectively, for $\lambda_1 > t$.
\end{Prop}

\noindent
\textbf{Proof.}
Follows from Proposition 4.
\
$\Box$

The central moments ($\mu_r$) and cumulants ($\kappa_r$)
of $X$ can be calculated from
\begin{eqnarray*}
\displaystyle
\mu_{r}=\sum_{k=0}^{r}(-1)^k \binom{r}{k} \mu_{1}^{\prime k} \mu_{r-k}^{\prime}
\qquad
\text{and}
\qquad
\displaystyle
\kappa_{r}=\mu_{r}^{\prime}-\sum_{k=1}^{r-1}\binom{r-1}{k-1} \kappa_{k} \mu_{r-k}^{\prime},
\end{eqnarray*}
respectively, where $\kappa_{1}=\mu^{\prime}_1$.
Note that $\kappa_{2}=\mu^{\prime}_2-\mu^{\prime 2}_1$,
$\kappa_{3}=\mu^{\prime}_3-3\mu^{\prime}_2 \mu^{\prime}_1+2\mu^{\prime 3}_1$,
$\kappa_4=\mu_4^{\prime}-4{\mu}_3^{\prime}\mu_1^{\prime}-3{\mu}_2^{\prime2}+12{\mu}_2^{\prime}{\mu}_1^{\prime 2}-6{\mu}_1^{\prime 4}$, etc.
The skewness $\gamma_1=\kappa_3/\kappa_2^{3/2}$ and kurtosis $\gamma_2=\kappa_4/\kappa_2^{2}$
follow from the second, third and fourth cumulants.

\begin{figure}[!htbp]
(a) \hspace{9cm} (b)

\includegraphics[width=0.45\textwidth]{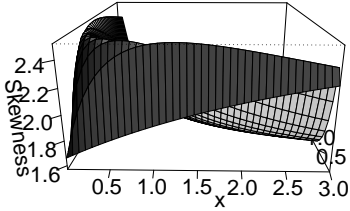}
\includegraphics[width=0.45\textwidth]{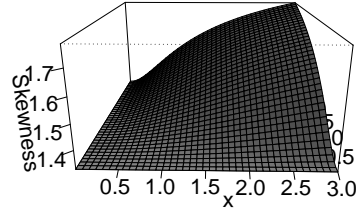}
\caption{Skewness for different values of $\lambda_2$ and $\alpha$ when  $\lambda_1=5$, $\theta=0.99$.
Skewness for different values of $\lambda_2$ and $\alpha$ when (a) $\lambda_1=5$, $\theta=0.99$,
(b) $\lambda_1=5$, $\theta=-0.99$.}
\label{skew1}
\end{figure}

\begin{figure}[!htbp]
(a)  \hspace{9cm}   (b)

\includegraphics[width=0.45\textwidth]{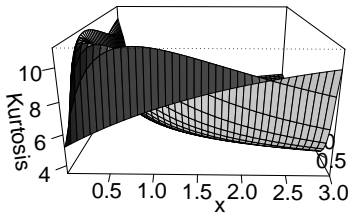}
\includegraphics[width=0.45\textwidth]{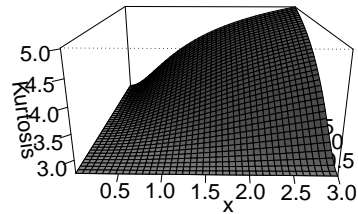}
\caption{Kurtosis for different values of $\lambda_2$ and $\alpha$ when  $\lambda_1=5$, $\theta=0.99$.}
\label{kurt1}
\end{figure}

Figures \ref{skew1} and \ref{kurt1} plot skewness and kurtosis
as functions of $\lambda_2$ and $\alpha$ when $\lambda_1$ and $\theta$ are fixed.
Skewness assumes values greater than $2$, implying that the GWED possesses  a wider range of skewness values
than the two-parameter WED due to Shakhatreh (2012).

The shape of many distributions can be usefully described by
conditional moments.
These  moments play an important role in measuring
inequality, for example, income quantiles, Lorenz curve and Bonferroni curve.
Proposition 6 derives the $n$th conditional moment and the conditional moment generating function of $X$.

\begin{Prop}
The $n$th conditional moment and the conditional moment generating function of $X$ can be expressed as
\begin{eqnarray*}
\displaystyle
m_n \left(x_0\right)
&=&
\displaystyle
p_{1} \frac {\displaystyle \gamma\left(n+1,\lambda_1 x_0\right)}{\displaystyle \lambda_1^n}+
p_{2} \frac {\displaystyle \gamma\left(n+1, \left(\lambda_1+\lambda_2\alpha\right) x_0\right)}{\displaystyle \left(\lambda_1+\lambda_2\alpha\right)^n}+
p_3 \frac {\displaystyle \gamma\left(n+1, \left( 2\lambda_1+\lambda_2\alpha \right) x_0\right)}{\displaystyle \left( 2\lambda_1+\lambda_2\alpha \right)^n}
\\
&&
\displaystyle
+p_4 \frac {\displaystyle \gamma\left(n+1, \left(\lambda_1+2\lambda_2\alpha\right) x_0\right)}{\displaystyle \left(\lambda_1+2\lambda_2\alpha\right)^n}+
p_5 \frac {\displaystyle \gamma\left(n+1,2 \left(\lambda_1+\lambda_2\alpha\right) x_0\right)}{\displaystyle 2^n \left(\lambda_1+\lambda_2\alpha\right)^n}
\end{eqnarray*}
and
\begin{eqnarray*}
\displaystyle
\mathcal{M^{\star}}(t)
&=&
\displaystyle
\frac {\displaystyle p_{1}\lambda_{1} \left[ 1 - e^{-\left(\lambda_1-t\right)x_0} \right]}{\displaystyle \lambda_{1} - t} +
\frac {\displaystyle p_{2}\left(\lambda_{1}+\lambda_{2}\alpha\right)
\left[ 1 - e^{-\left(\lambda_{1}+\lambda_{2}\alpha -t\right)x_0}\right]}{\displaystyle \lambda_{1}+\lambda_{2}\alpha - t}
\\
&&
\displaystyle
+\frac {\displaystyle p_{3} \left( 2\lambda_{1}+\lambda_{2}\alpha \right)
\left[ 1 - e^{-\left( 2\lambda_{1}+\lambda_{2}\alpha - t\right)x_0} \right]}{\displaystyle 2\lambda_{1}+\lambda_{2}\alpha -t} +
\frac {\displaystyle p_{4} \left( \lambda_{1}+2\lambda_{2}\alpha \right)
\left[ 1 - e^{-\left( \lambda_{1}+2\lambda_{2}\alpha -t \right)x_0} \right]}
{\displaystyle \lambda_{1}+2\lambda_{2}\alpha -t}
\\
&&
\displaystyle
+\frac {\displaystyle 2p_{5} \left( \lambda_{1}+\lambda_{2}\alpha \right)
\left[ 1 - e^{-\left( 2\lambda_{1}+2\lambda_{2}\alpha -t \right)x_0} \right]}{\displaystyle  2\lambda_{1}+2\lambda_{2}\alpha -t},
\end{eqnarray*}
respectively, for $\lambda_1 > t$,
where $\gamma(a,x) = \displaystyle \int_0^x t^{a-1} e^{-t} {\rm d}t$ denotes the lower incomplete gamma function.
\end{Prop}

\noindent
\textbf{Proof.}
Follows from Proposition 4.
\
$\Box$

Probability weighted moments (PWMs) formally defined as
\begin{eqnarray*}
\displaystyle
\tau_{s,r} = E\left\{X^s \left[F^{w}\left(X\right)\right]^r\right\} = \int_0^{+\infty}x^s F^{w}(x)^r f^{w}(x) {\rm d}x
\end{eqnarray*}
are used to summarize distributions.
They are also used for estimation of parameters
especially when the inverse cdf cannot be expressed explicitly.
Proposition 7 derives an expression for the PWMs of $X$.

\begin{Prop}
The PWMs of $X$ can be expressed as
\begin{eqnarray*}
\displaystyle
\tau_{s,r}
&=&
\displaystyle
\sum_{i\ge0}c_{r,i} (s+r+i)! \Bigg[\frac {\displaystyle p_1}{\displaystyle \lambda_1^{s+r+i}}+
\frac {\displaystyle p_{2}}{\displaystyle \left(\lambda_1+\lambda_2\alpha\right)^{s+r+i}}+
\frac {\displaystyle p_3}{\displaystyle \left(2\lambda_1+\lambda_2\alpha\right)^{s+r+i}}
\\
&&
\displaystyle
\qquad \qquad \qquad \qquad \qquad \qquad
+\frac {\displaystyle p_4}{\displaystyle \left(\lambda_1+2\lambda_2\alpha\right)^{s+r+i}}+
\frac {\displaystyle p_5}{\displaystyle 2^{s+r+i} \left(\lambda_1+\lambda_2\alpha\right)^{s+r+i}}\Bigg],
\end{eqnarray*}
where $c_{r, 0}=d_0^r$ and
$c_{r, m}=\left(a_0 m\right)^{-1}\displaystyle \sum_{k=1}^m \left[k(r+1)-m\right] d_k c_{r, m-k}$ for $m\ge1$.
\end{Prop}

\noindent
\textbf{Proof.}
By definition,
\begin{eqnarray*}
\displaystyle
\tau_{s,r} = K\int_0^{+\infty}x^{s+r}\left[\sum_{i\ge0}d_i x^i\right]^r
f^{w} \left(x;\lambda_{1}, \lambda_{2}, \alpha, \theta\right){\rm d}x.
\end{eqnarray*}
By equation (0.314)  in Gradshteyn and Ryzhik (2000),
\begin{eqnarray*}
\displaystyle
\left[\sum_{i\ge0}d_i x^i\right]^r=\sum_{i\ge0}c_{r,i} x^i,
\end{eqnarray*}
where $c_{r, 0} = d_0^r$ and $c_{r, m}=\left(a_0 m\right)^{-1} \displaystyle\sum_{k=1}^m \left[k(r+1)-m\right] d_k c_{r,m-k}$ for $m\ge 1$.
The result follows by using \eqref{pdfGWEmix}.
\
$\Box$

We now invert \eqref{P1} to obtain a power series expansion for the
quantile function say $x=Q(u)$ of the GWED.
We shall use the Lagrange theorem.
We assume that the power series expansion
\begin{eqnarray*}
\displaystyle
g=F^{w}(x)=g_0+\sum_{n=1}^{+\infty}f_n \left(x-x_0\right)^n
\end{eqnarray*}
holds, where $F^{w}(x)$ is analytic at a point $x_0$.
Then, the inverse function $x=Q(u)=\left(F^{w}(u)\right)^{-1}$ exists in the neighborhood of some point $u=u_0$.

\begin{Prop}
The quantile function of $X$ can be expressed as
\begin{eqnarray}
\label{qq}
\displaystyle
Q(u)=\sum_{n=1}^{+\infty}b_n u^n,
\end{eqnarray}
where $b_n= q_{n, n-1} / \left(n d_0^{n} \right)$,
$q_{n, i}=i^{-1} \displaystyle \sum_{m=1}^{i} \left[ m (n+1) - i \right] p_m q_{n, i-m}$,
$q_{n, 0}=p_0^n=1$ and $p_i =-m_0^{-1} \displaystyle \sum_{k=1}^i d_{k+1} p_{i-k}$, $i\geq 1$ with $p_0 = 1$.
\end{Prop}

\noindent
\textbf{Proof.}
According to Markushevich (1965, volume 2, page 88),  a power series of $x=Q(u)$ is
\begin{eqnarray*}
\displaystyle
x=Q(u)=x_0+\sum_{n=1}^{+\infty} h_n \left(u-u_0\right)^n,
\end{eqnarray*}
where
\begin{eqnarray*}
\displaystyle
h_n = \frac {\displaystyle 1}{\displaystyle n!}
\frac {\displaystyle d^{n-1}}{\displaystyle d z^{n-1}}\left\{ \left[\psi(x)\right]^{n}\right\}
\bigg{|}_{x=x_0}
\quad
\mbox{and}
\quad
\displaystyle
\psi(x)=\frac {\displaystyle x-x_0}{\displaystyle F^{w}(x)-F^{w}\left(x_0\right)}.
\end{eqnarray*}
Setting $x_0=0$ and $u_0=0$, and using \eqref{P1}, we have
\begin{eqnarray*}
\displaystyle
\psi(x) = \frac {\displaystyle x}{\displaystyle F^{w}(x)-0} = \frac {\displaystyle 1}{\displaystyle \sum\limits_{i\ge0}d_i x^i}.
\end{eqnarray*}
The inverse of the power series $\displaystyle \sum_{i\ge0} d_i x^i$ follows from Gradshteyn and Ryzhik (2000, equation (0.313)):
\begin{eqnarray*}
\displaystyle
\psi(x)=\frac {\displaystyle 1}{\displaystyle \sum\limits_{i\ge0} d_i x^i}=\frac {\displaystyle 1}{\displaystyle d_0}\sum_{i\ge0} p_i x^{i},
\end{eqnarray*}
where the coefficients $p_i$ can be determined from
$p_i =-m_0^{-1} \displaystyle\sum_{k=1}^i d_{k+1} p_{i-k}$, $i\geq 1$ with $p_0= 1$.
Using Gradshteyn and Ryzhik (2000, equation (0.314)),
we can write $\psi(x)^{n} = \frac {1}{d_0^n} \displaystyle\sum\limits_{i\ge0} q_{n,i} x^i$,
where $q_{n,i}$ for $i=1,2,\ldots$ are given by
$q_{n,i}=i^{-1} \displaystyle \sum_{m=1}^{i} \left[ m (n+1)-i \right] p_m q_{n,i-m}$,
and $q_{n,0}=p_0^n=1$.
The $q_{n,i}$ can be determined from $q_{n,0}, \ldots, q_{n,i-1}$ and therefore from $p_0, \ldots, p_i$.
The derivative of order $(n-1)$ of  $\psi(x)^{n}$ can be expressed as
\begin{eqnarray*}
\displaystyle
h_n = \frac {\displaystyle 1}{\displaystyle n!}
\frac {\displaystyle d^{n-1}}{\displaystyle d x^{n-1}}
\left\{ \left[\psi(x)\right]^{n}\right\}\bigg{|}_{x=0} = \frac {\displaystyle q_{n,n-1}}{\displaystyle n d_0^{n}}.
\end{eqnarray*}
Hence, the power series for the quantile function reduces to \eqref{qq}.
\
$\Box$

\subsection{Maximum likelihood estimation}

Here, we consider estimation of the parameters $\boldsymbol{\eta} = \left(\lambda_{1}, \lambda_{2}, \alpha, \theta \right)$
of the GWED by the method of maximum likelihood.
We suppose $x_1,x_2,\ldots,x_n$ is a random sample from the GWED.
Then, the log-likelihood function is
\begin{eqnarray*}
\displaystyle
\ell \left(\boldsymbol{\eta}\right)
&=&
\displaystyle
n\left\{\log \lambda_1-\log\left[\frac {\displaystyle \lambda_2\alpha}{\displaystyle \lambda_1+\lambda_2\alpha} +
\theta\lambda_1\left(\frac {\displaystyle 2}{\displaystyle 2\lambda_1+\lambda_2\alpha} +
\frac {\displaystyle 1}{\displaystyle \lambda_1+2\lambda_2\alpha} -
\frac {\displaystyle 2}{\displaystyle \lambda_1+\lambda_2\alpha}\right)\right]\right\}
\\
&&
\qquad \qquad
\displaystyle
-\lambda_1\sum_{i=1}^n x_i+
\sum_{i=1}^n \log\left(1 - e^{-\lambda_2 \alpha x_i}\right) +
\sum_{i=1}^n \log\left[1-\theta  e^{-\lambda_2 \alpha x_i}
\left(1-2e^{-\lambda_1 x_i}\right)\right].
\end{eqnarray*}
Differentiating $\ell \left(\boldsymbol{\eta}\right)$ with respect to $\lambda_{1}$, $\lambda_{2}$, $\alpha$ and $\theta$,
we obtain the normal equations
\begin{eqnarray*}
&&
\displaystyle
\frac {\displaystyle \partial \ell \left(\boldsymbol{\eta}\right)}{\displaystyle \partial \lambda_1}=
\frac {\displaystyle n}{\displaystyle \lambda_1} - n \lambda_2\alpha \frac {\displaystyle A}{\displaystyle B} - \sum_{i=1}^n x_i -
2\theta\sum_{i=1}^n x_i \frac {\displaystyle e^{-\left(\lambda_1+\lambda_2\alpha\right)x_i}}
{\displaystyle 1-\theta e^{-\lambda_2\alpha x_i} \left( 1 - 2e^{-\lambda_1 x_i} \right)}=0,
\\
&&
\displaystyle
\frac {\displaystyle \partial \ell\left(\boldsymbol{\eta}\right)}{\displaystyle \partial \lambda_2} =
n \lambda_1\alpha \frac {\displaystyle A}{\displaystyle B} +
\alpha\sum_{i=1}^n \frac {\displaystyle x_ie^{-\lambda_2\alpha x_i}}{\displaystyle 1-e^{-\lambda_2 \alpha x_i}} +
\theta\alpha \sum_{i=1}^n x_i \frac {\displaystyle e^{-\lambda_2\alpha x_i} - 2e^{-\left( \lambda_1+\lambda_2 \right)x_i}}
{\displaystyle 1-\theta e^{-\lambda_2\alpha x_i} \left(1-2e^{-\lambda_1 x_i} \right)} = 0,
\\
&&
\displaystyle
\frac {\displaystyle \partial \ell\left(\boldsymbol{\eta}\right)}{\displaystyle \partial \alpha} =
-n \lambda_1\lambda_2 \frac {\displaystyle A}{\displaystyle B} -
\lambda_2\sum_{i=1}^n \frac {\displaystyle x_i e^{-\lambda_2\alpha x_i}}{\displaystyle 1-e^{-\lambda_2 \alpha x_i}} +
\theta\lambda_2 \sum_{i=1}^n x_i \frac {\displaystyle e^{-\lambda_2\alpha x_i} -
2 e^{-\left( \lambda_1+\lambda_2 \right)x_i}}{\displaystyle 1 - \theta e^{-\lambda_2\alpha x_i}\left( 1 - 2e^{-\lambda_1x_i} \right)}=0,
\\
&&
\displaystyle
\frac {\displaystyle \partial \ell \left(\boldsymbol{\eta}\right)}{\displaystyle \partial \theta} =
-\frac {\displaystyle n\lambda_1 \left(\frac {2}{2\lambda_1+\lambda_2\alpha} + \frac {1}{\lambda_1+2\lambda_2\alpha} -
\frac {2}{\lambda_1+\lambda_2\alpha}\right)}{\displaystyle B} -
\sum_{i=1}^n \frac {\displaystyle e^{-\lambda_2 \alpha x_i}\left(1-2e^{-\lambda_1 x_i}\right)}
{\displaystyle 1 - \theta e^{-\lambda_2 \alpha x_i}\left(1-2e^{-\lambda_1 x_i}\right)} = 0,
\end{eqnarray*}
where
\begin{eqnarray*}
\displaystyle
A = -\frac {\displaystyle 1}{\displaystyle \left( \lambda_1+\lambda_2\alpha \right)^2} +
2\theta\left[ \frac {\displaystyle 1}{\displaystyle \left( 2\lambda_1+ \lambda_2\alpha\right)^2} +
\frac {\displaystyle 1}{\displaystyle \left(\lambda_1+2\lambda_2\alpha\right)^2} -
\frac {\displaystyle 1}{\displaystyle \left( \lambda_1+\lambda_2\alpha \right)^2}\right]
\end{eqnarray*}
and
\begin{eqnarray*}
\displaystyle
B=\frac {\displaystyle \lambda_2\alpha}{\displaystyle \lambda_1+\lambda_2\alpha}+
\theta\lambda_1\left(\frac {\displaystyle 2}{\displaystyle 2\lambda_1+ \lambda_2\alpha} +
\frac {\displaystyle 1}{\displaystyle \lambda_1+2\lambda_2\alpha} -
\frac {\displaystyle 2}{\displaystyle \lambda_1+\lambda_2\alpha}\right).
\end{eqnarray*}
The maximum likelihood estimates say $\boldsymbol{\widehat{\eta}} = \left( \widehat{\lambda}_{1}, \widehat{\lambda}_{2},
\widehat{\alpha}, \widehat{\theta} \right)$
are the simultaneous solutions of the normal equations.
These equations do not yield explicit solutions.
Hence, the maximum likelihood estimates must be obtained numerically.

According to Cox and Hinkley (1979), the asymptotic distribution of
$\boldsymbol{\widehat{\eta}}$ can be approximated by the
multivariate normal distribution, $N_{4}\left( \textbf{0},\left[\textbf{J} \left(\boldsymbol{\widehat{\eta}} \right)\right]^{-1} \right)$,
where $\textbf{J} \left(\boldsymbol{\widehat{\eta}}\right)$ denotes
the inverse of the observed information matrix evaluated at $\boldsymbol{\widehat{\eta}}$.
Due to algebraic complexity and in order to save space,
we have not reported the expression for $\textbf{J}\left( \boldsymbol{\widehat{\eta}} \right)$.
This approximation can be used to construct approximate confidence intervals and hypothesis
tests for $\lambda_{1}$, $\lambda_{2}$, $\alpha$ and $\theta$.

Numerical calculations not reported here showed that
the surface of the $\ell \left(\boldsymbol{\eta}\right)$ was smooth.
Numerical routines for maximization of  $\ell \left(\boldsymbol{\eta}\right)$
were able to locate the maximum for a wide range of starting values.
However, to ease computations it is useful to have reasonable starting values.
These can be
obtained, for example, by the method of moments.
For $r=1,2,3,4$ let $m_r = \displaystyle \frac {1}{n}\sum_{i=1}^n x_i^r$ denote the first four sample moments.
Equating these moments with the theoretical versions given in Section 3.1,
we have $m_r = E \left(X^r\right)$ for $r = 1, 2, 3, 4$.
These equations can be solved simultaneously
to obtain the moments estimates.

\subsection{Application}

Here, we illustrate the flexibility of the GWED using two biochemical real  data sets:
C-reactive protein (CRP) data and insulin data.
These data were collected
through  a study whose aim was to examine the correlation of the RBP4 (retinol-binding protein-4)
with anthropometric measurements (BMI-body mass index and WC-waist circumference),
insulin resistance, metabolic and kidney parameters and inflammation.
The participants
of the study were 128 obese diabetic patients.
This research was carried out in Primary Health Care Center (Center of Laboratory Diagnostics)
and  Center of Clinical-Laboratory Diagnostics, Clinical Center of Montenegro.
These data sets are original in that they have not been analyzed in the statistics literature before.

We fitted the GWED to both data sets by the method of maximum likelihood.
The fit of the GWED was compared with that of the WED due to Gupta and Kundu (2009).
As criteria for comparison, we used the Akaike information
criterion (AIC) and the $p$-value of the Kolmogorov-Smirnov test.

Both distributions were fitted by
executing the R function {\sf optimx} (Nash and  Varadhan, 2011) for a wide range of starting values.
This sometimes resulted in more than one maximum, but at least one maximum was
identified each time.
In cases of more than one maximum, we took the maximum likelihood estimates to correspond to the largest of the maxima.

\subsubsection{Application 1: CRP data}

CRP is a protein found in the blood, the levels of which rise in response to inflammation
(i.e., CRP is an acute-phase protein).
The CRP gene is located on the first chromosome (1q21-q23).
Table \ref{D1} gives descriptive statistics of the CRP data set.
We see that the data has positive skewness and kurtosis greater than that of the normal distribution.

\begin{table}
\caption{Descriptive statistics of the CRP data.}
\label{D1}
\begin{tabular}{lllllll}
\hline
\noalign{\smallskip}
Minimum & Maximum & Mean & Median & SD &  Skewness & Kurtosis
\\
\noalign{\smallskip}
\hline
\noalign{\smallskip}
0.15 & 7.85 & 1.8867 & 1.485 & 1.5450 &  1.1777& 4.1425
\\
\noalign{\smallskip}
\hline
\end{tabular}
\end{table}

The MLEs of the model parameters, their standard errors, the AIC values,
and the $p$-values of the Kolmogorov-Smirnov test for the CRP data are reported in Table \ref{D2}.
Based on the AIC values and the $p$-values, we see that the GWED provides a better fit than the WED for the CRP data.

\begin{table}
\caption{MLEs of the model parameters, AIC values and $p$-values  for the CRP data.}
\label{D2}
\begin{tabular}{llllllll}
\hline
\noalign{\smallskip}
Model & $\lambda_1$ & $\lambda_2$ & $\alpha$ & $\theta$ &  $\lambda$ & AIC & $p$-value
\\
\noalign{\smallskip}
\hline
\noalign{\smallskip}
GWED & 0.7280 & 0.8316  & 0.8121& 0.7991&  -- & 574.84 & 0.819
\\
& (0.0874)& (0.0792) & (0.1010)&(0.0936)&&&
\\
WED & -- & -- & 8.5271 & -- &  0.5857 & 777.76 & 0
\\
& & &(1.0024) & & (0.08715) & &
\\
\noalign{\smallskip}
\hline
\end{tabular}
\end{table}

\begin{figure}[!htbp]
(a)  \hspace{9cm}  (b)

\includegraphics[width=0.45\textwidth]{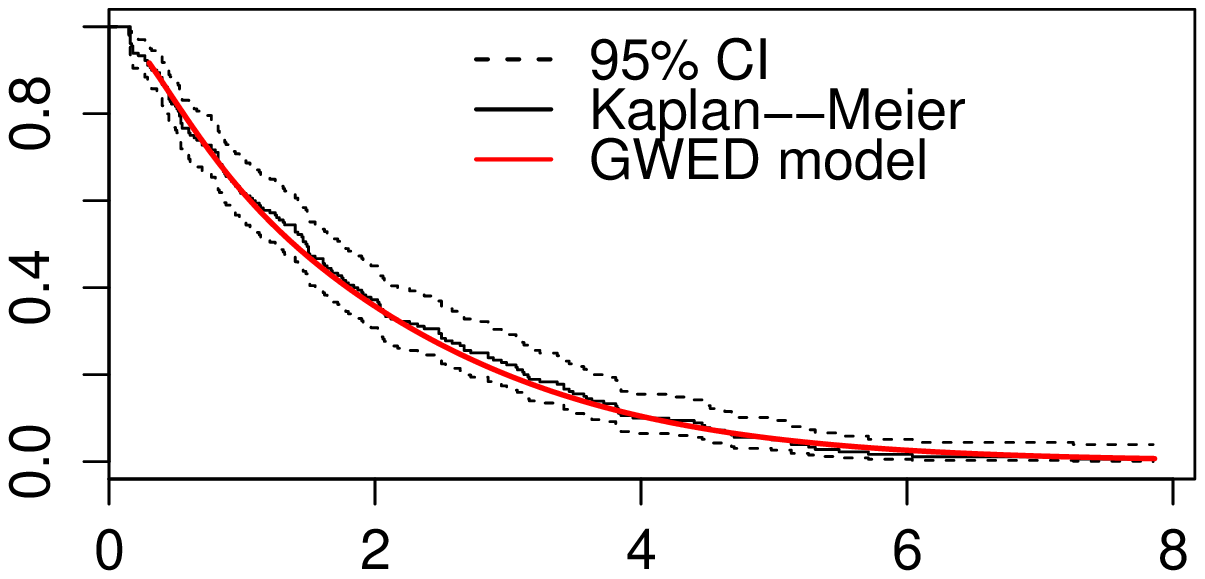}
\includegraphics[width=0.45\textwidth]{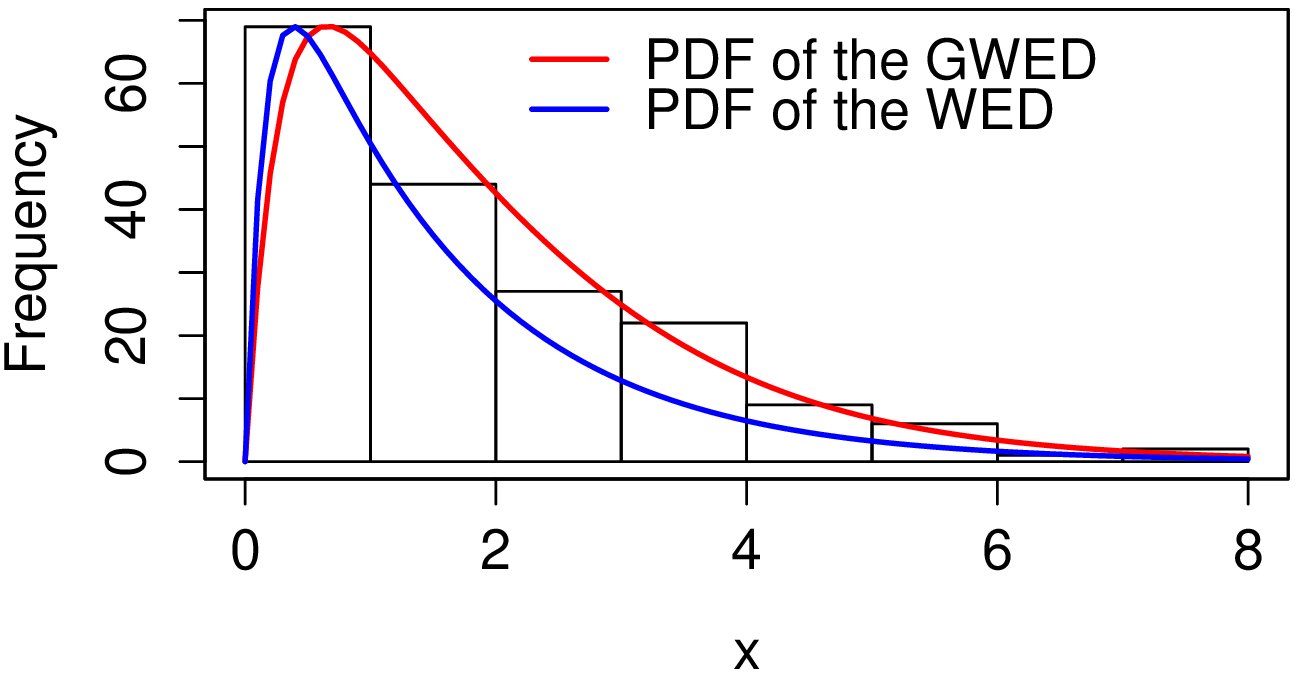}
\caption{(a) Fitted  GWE survival function and the empirical survival function for the CRP data;
(b) Fitted GWE and WE pdfs and the empirical pdf for the CRP data.}
\label{KM1}
\end{figure}

\begin{figure}[!htbp]
(a) \hspace{9cm}  (b)

\includegraphics[width=0.45\textwidth]{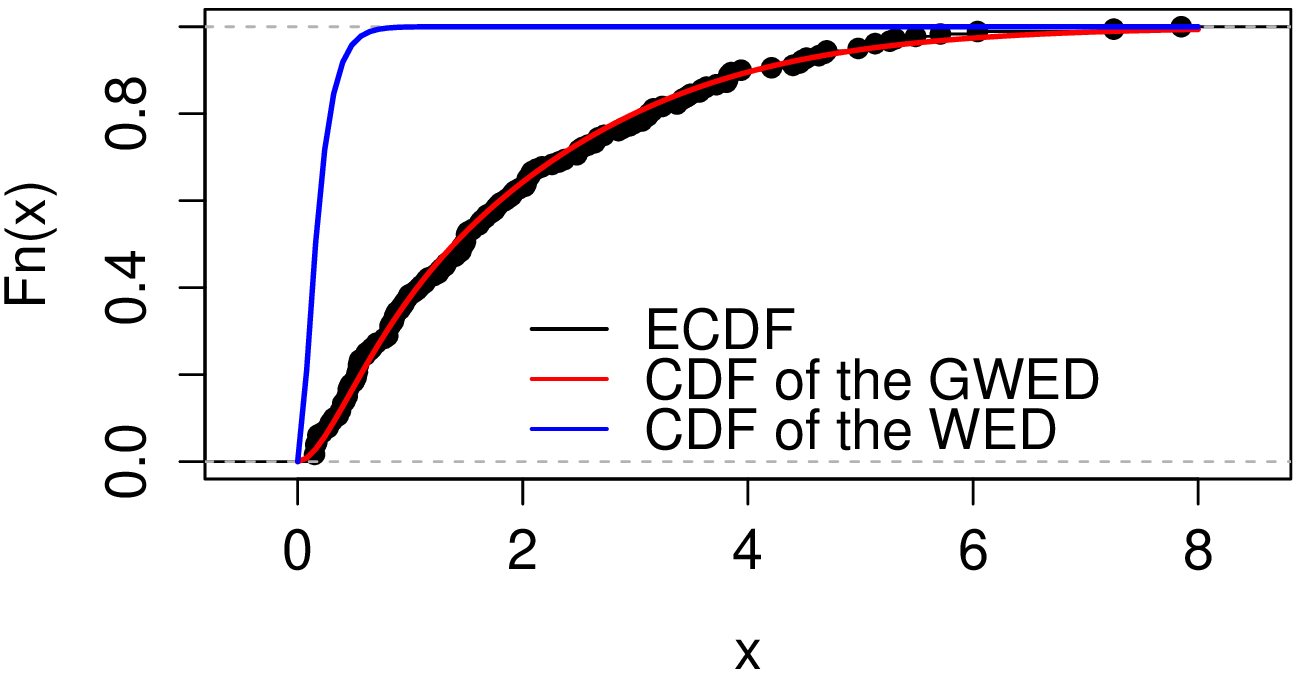}
\includegraphics[width=0.45\textwidth]{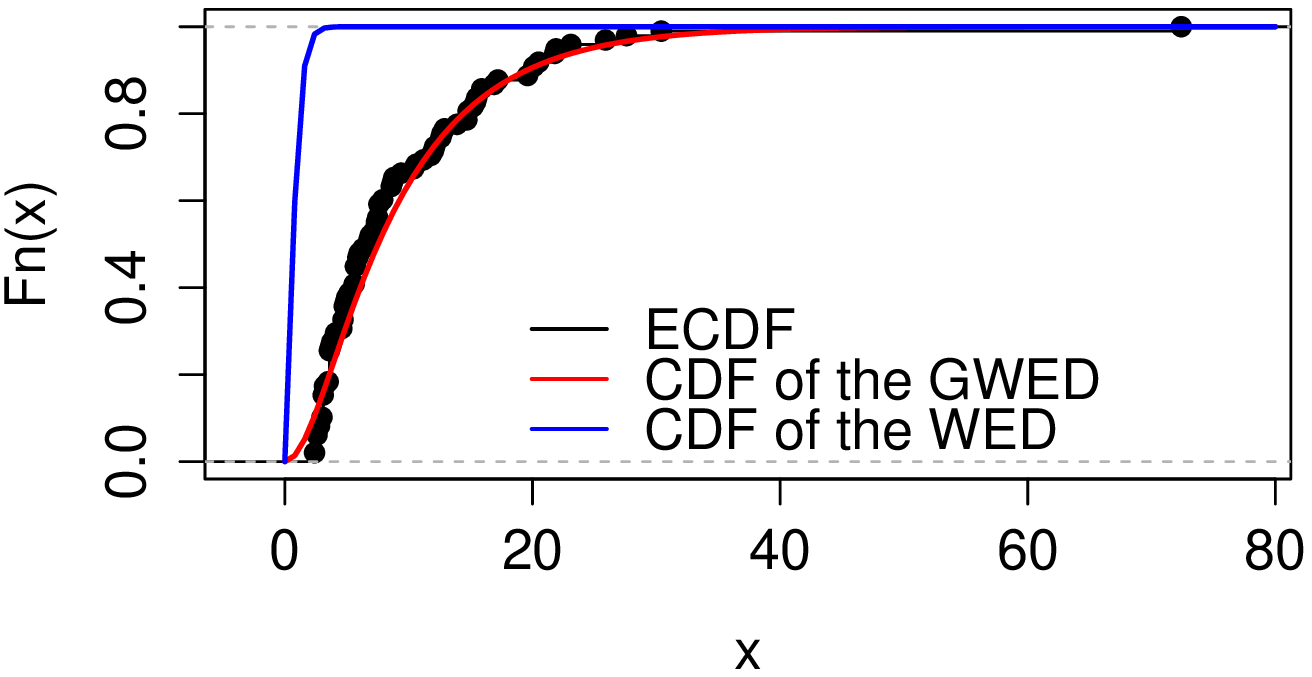}
\caption{(a) Fitted GWE and WE cdfs and the empirical cdf for the  CRP data;
(b) Fitted GWE and WE cdfs and the empirical cdf for the insulin data.}
\label{ecdf}
\end{figure}

Figure \ref{KM1}(a)  plots the survival function for the fitted GWED and the empirical survival function for the CRP data.
Figure \ref{KM1}(b) plots the pdfs for the fitted GWED and WED and the empirical pdf for the CRP data.
Figure \ref{ecdf}(a) plots the cdfs for the fitted GWED and WED and the empirical cdf for the CRP data.
All these figures suggest that the GWED provides a better fit to the CRP data.

\subsubsection{Application 2: Insulin data}

Insulin is a peptide hormone, produced by beta cells of the pancreas,
and is central to regulating carbohydrate and fat metabolism in the body.
Insulin causes cells in the liver, skeletal muscles, and fat tissue to absorb glucose from the blood.
In the liver and skeletal muscles, glucose is stored as glycogen, and in fat cells (adipocytes) it is stored as triglycerides.
The descriptive statistics of the insulin data are given in Table \ref{D3}.
This data is also positively skewed and has kurtosis greater than that of the normal distribution.

\begin{table}
\caption{Descriptive statistics of the insulin data.}
\label{D3}
\begin{tabular}{lllllll}
\hline
\noalign{\smallskip}
Minimum & Maximum & Mean & Median & SD &  Skewness & Kurtosis
\\
\noalign{\smallskip}
\hline
\noalign{\smallskip}
2.4 & 72.4 & 9.616 & 6.8 & 9.1306 &  1.7050 & 24.158
\\
\noalign{\smallskip}
\hline
\end{tabular}
\end{table}

The MLEs of the model parameters, their standard errors, the AIC values,
and the $p$-values of the Kolmogorov-Smirnov test for the insulin data are reported in Table \ref{D3}.
Based on the AIC values and the $p$-values, we see again that the GWED provides a better fit than the WED for the insulin data.

\begin{table}
\caption{MLEs of the model parameters, AIC values and $p$-values  for the insulin data.}
\label{D4}
\begin{tabular}{llllllll}
\hline
\noalign{\smallskip}
Model & $\lambda_1$ & $\lambda_2$ & $\alpha$ & $\theta$ &  $\lambda$ & AIC & $p$-value
\\
\noalign{\smallskip}
\hline
\noalign{\smallskip}
GWED & 0.1352 & 0.8885 & 0.4019 & -0.3499 &  -- & 617.56 & 0.2225
\\
& (0.0054)&(0.1272)&(0.0157)&(0.0479)&&&
\\
WED & -- & -- & 2.3631 & -- &  0.1349 & 933.01 & 0
\\
& & &(0.9987)& &(0.0241)& &
\\
\noalign{\smallskip}
\hline
\end{tabular}
\end{table}

\begin{figure}[!htbp]
(a) \hspace{9cm}  (b)

\includegraphics[width=0.45\textwidth]{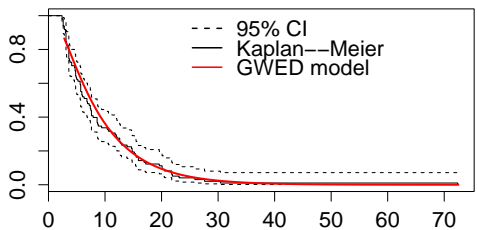}
\includegraphics[width=0.45\textwidth]{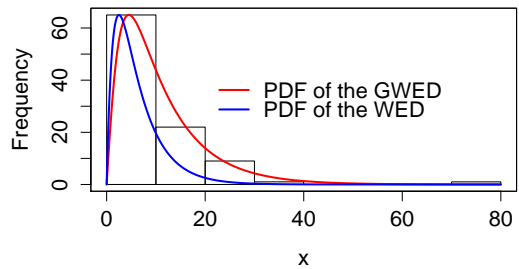}
\caption{(a) Fitted  GWE survival function and the empirical survival function for the insulin data;
(b) Fitted GWE and WE pdfs and the empirical pdf for the insulin data.}
\label{KM2}
\end{figure}

Figure \ref{KM2}(a)  plots the survival function for the fitted GWED and the empirical survival function for the insulin data.
Figure \ref{KM2}(b) plots the pdfs for the fitted GWED and WED and the empirical pdf for the insulin data.
Figure \ref{ecdf}(b) plots the cdfs for the fitted GWED and WED and the empirical cdf for the insulin data.
All these figures suggest that the GWED again provides a better fit to the insulin data.

\section{Concluding remarks and future research}

In this paper, we have proposed an extension of Azzalini's method.
As an illustration of this extension, we have introduced
a four-parameter  Generalized Weighted Exponential distribution (GWED).
This new  distribution  generalizes of the weighted exponential distribution due to Gupta and Kundu (2009).
The method that we have proposed can be used to generalize any distribution.

We have studied mathematical properties of the GWED.
We have shown that the GWED can be interpreted as a hidden truncated distribution (Arnold and Beaver, 2000a).
The dependence parameter rather than the marginal parameters
plays an important role in making the GWED flexible in terms of skewness,
kurtosis, shape of the hazard rate function and other characteristics.

We have discussed maximum likelihood estimation
of the parameters of the GWED and
provided two real data applications.
The applications show that the GWED can be used quite
effectively to give better fits than the WED.
We hope that the GWED may attract wider applications in statistics.

In conclusion, it is important to highlight that although
this paper has mainly concentrated on the GWED,
our method can be used for any distribution.
For instance, one can consider a general class of cdfs defined by
\begin{eqnarray*}
\displaystyle
F\left( x\right) =F\left( x;a,b,c\right) =\left[ a\xi \left( x\right) +b\right]^{c},
\end{eqnarray*}
where $a\neq 0$, $b \neq 0$ and $c\neq 0$ are parameters, and
$\xi \left( x\right) $  is a monotonic and differentiable function on $\left( F^{-1}\left( 0\right), F^{-1}\left( 1\right) \right)$.
This class for proper choices of $a$, $b$, $c$ and $\xi \left( \cdot \right)$
can contain various distributions studied in the literature
like the power function, Dagum (Burr type III), Pareto, inverse Weibull, reflected
exponential and rectangular distributions.

\section*{Acknowledgement}

We are grateful to Dr. Aleksandra Klisi\'c, a biochemist in Primary Health Unit in
Podgorica, Montenegro, for allowing us to use the CRP and insulin data.

We would like to thank to the Editor and to anonymous referees whose comments
greatly improved quality of the paper.

\begin{center}
\textbf{References}
\end{center}

\noindent
Arnold, B.C. and  Beaver, R.J. (2000a).
Hidden truncation models.
Sankhy\=a, 62, pp. 23-35.

\noindent
Arnold, B.C. and  Beaver, R.J. (2000b).
The skew-Cauchy distribution.
Statistics and Probability Letters, 49, pp. 285-290.

\noindent
Aryal, G. and Nadarajah, S. (2005).
On the skew-Laplace distribution.
Journal of Information and Optimization Sciences, 26, pp. 205-217.

\noindent
Azzalini, A. (1985).
A class of distributions which includes the normal ones.
Scandinavian Journal of Statistics, 12, pp. 171-178.

\noindent
Cox, D.R. and Hinley, D.V. (1979).
Theoretical Statistics.
Chapman and Hall, London.

\noindent
Domma, F. and Giordano, S. (2013).
A copula-based approach to account for dependence in stress-strength models.
Statistical Papers, 54, pp. 807-826.

\noindent
Gradshteyn, I.S. and Ryzhik, I.M. (2000).
Table of Integrals, Series, and Products.
Academic Press, San Diego.

\noindent
Gupta, A.K., Chang, F.C. and Haung, W.J. (2002).
Some skew-symmetric models.
Random Operators and Stochastic Equations, 10, pp. 133-140.

\noindent
Gupta, R.D. and Kundu, D. (2009).
A new class of weighted exponential distributions.
Statistics, 43, pp. 621-634.

\noindent
Kotz, S., Lumelskii, Y. and Pensky, M. (2003).
The Stress-Strength Model and Its Generalizations: Theory and Applications.
World Scientific Publishing, Singapore.

\noindent
Mahdy, M. (2011).
A class of weighted gamma distributions and its properties.
Economic Quality Control, 26, pp. 133-144.

\noindent
Mahdy, M. (2013).
A class of weighted Weibull distributions and its properties.
Studies in Mathematical Sciences, 6, pp. 35-45.

\noindent
Markushevich, A.I. (1965).
Theory of Functions of a Complex Variable.
Chelsea Publication Company.

\noindent
Morgenstern, D. (1956).
Einfache Beispiele zweidimensionaler Verteilungen.
Mitteilingsblatt f\"{u}r Mathematische Statistik, 8, pp. 234-235.

\noindent
Nadarajah, S. (2009).
On the skew-logistic distribution.
AStA Advances in Statistical Analysis, 93, pp. 187-203.

\noindent
Nadarajah, S. (2014).
Expansions for bivariate copulas.
Submitted.

\noindent
Nash, J.C. and  Varadhan, R. (2011).
Unifying optimization algorithms to aid software system  users: {\sf optimx} for R.
Journal of Statistical Software, 43, pp. 1-14.

\noindent
Shakhatreh, M.K. (2012).
A two-parameter of weighted exponential distributions.
Statistics and Probability Letters, 82, pp. 252-261.

\noindent
Titterington, D.M.,  Smith, A.F.M. and Makov, U.E. (1985).
Statistical Analysis of Finite Mixture Distributions.
John Wiley and Sons, New York.

\end{document}